\documentclass[12pt]{amsart}
\usepackage{amssymb,amsmath,amscd,graphicx,
latexsym,amsthm}
\usepackage{amssymb,latexsym,eufrak,amsmath,amscd,graphicx}
  \usepackage[all]{xy}
\theoremstyle{plain}
\newtheorem{theorem}{Theorem}[section]

\newtheorem{proposition}[theorem]{Proposition}

\newtheorem{lemma}[theorem]{Lemma}
\newtheorem{variant}[theorem]{Variant}

\newtheorem{theoremalpha}{Theorem}
\newtheorem{corollaryalpha}[theoremalpha]{Corollary}

\theoremstyle{definition}

\newtheorem{remark}[theorem]{Remark}
\newtheorem{example}[theorem]{Example}

\newcommand{\lra}{\longrightarrow}
\newcommand{\noi}{\noindent}

\newcommand{\CC}{\mathbf{C}}

\newcommand{\OO}{\mathcal{O}}

\newcommand{\II}{\mathcal{I}}

\newcommand{\fra}{\mathfrak{a}}
\newcommand{\frb}{\mathfrak{b}}

\newcommand{\bull}{_{\bullet}}

\newcommand{\frakm}{\mathfrak{m}}
\newcommand{\frq}{\frak{q}}

\newcommand{\eps}{\varepsilon}

\newcommand{\MI}[1]{\mathcal{J}  ( {#1}
) }

\newcommand{\MMI}[2]{\MI{ {#1} \, , \, {#2}}}

\newcommand{\ord}{\textnormal{ord}}

\newcommand{\pr}{\prime}

\newcommand{\tnparenonly}[1]{\textnormal{(}{#1}\textnormal{)}}

\newcommand{\MIJ}  {\mathcal{J}}
\newcommand{\Tor}{\textnormal{Tor}}
\newcommand{\Image}{\textnormal{Im}}

\begin{document}

\title{Local syzygies of multiplier ideals}

\author{Robert Lazarsfeld}
\address{Department of Mathematics, University of Michigan, Ann Arbor, MI
 48109}
\email{{\tt rlaz@umich.edu}}
\thanks{Research of the first author partially supported by NSF grant DMS
0139713}

\author{Kyungyong Lee}
\address{Department of Mathematics, University of Michigan, Ann Arbor, MI
 48109}
\email{{\tt kyungl@umich.edu}}

\maketitle

\section*{Introduction}
The purpose of this note is to prove an elementary but somewhat
surprising result concerning the syzygies of multiplier ideals. It follows
in particular that in dimensions three or higher, multiplier ideals are
very special among all integrally closed ideals.

Let $X$ be a smooth complex algebraic variety of dimension $ d$,
and let
$\frb \subseteq \OO_X$ be an ideal sheaf. Given a rational or real number
$c > 0$ one can construct the \textit{multiplier ideal} \[
\MI{\frb^c} \ = \MI{X, \frb^c } \ \subseteq \ \OO_X\]
of $\frb$ with weighting coefficient $c$. This is a new ideal on $X$ that
measures in a somewhat subtle manner the singularities of functions $f
\in \frb$. In recent years,  multiplier ideals  have
found many applications in local and global algebraic geometry (cf. \cite{Dem1},
\cite{AS},  \cite{EL}, \cite{Siu1}, \cite{ELS}, \cite{Siu2}, \cite{Tak}, \cite{HM1}, \cite{HM2}, \cite{Dem2}, \cite{PAG}). 

Because of their importance, there has been some interest in trying to
understand how general or special multiplier ideals may   be among all
ideal sheaves. Multiplier ideals are always integrally closed, but up to
now they have not been known to satisfy any other local properties. In
fact, Favre--Jonsson \cite{FJ} and  Lipman--Watanabe \cite{LW} proved that in
dimension $d = 2$, every integrally closed ideal can locally be realized
as a multiplier ideal. 

The first examples of integrally closed non-multiplier
ideals in dimension
$\ge 3$ were discovered by the second author, who recently gave some
quite delicate geometric arguments to show that the ideal of a suitable
number of general lines through the origin in $\CC^3$ couldn't arise as a
multiplier ideal. However the construction didn't pinpoint any general
features of multiplier ideals that might be violated: rather the idea
was to follow  a potential resolution of singularities of the data with
enough care that one could eventually get a
contradiction.\footnote{We remark that these examples are in fact not covered by
Theorem \ref{TheoremA}.}

 Our main result shows that multiplier ideals satisfy some possibly unexpected
properties of an algebraic nature.  In the following, we work in the local
ring
$(\OO,
\frakm)$ of  
$X$  at a point $x\in X$, and as above $d = \dim X$. 
\begin{theoremalpha} \label{TheoremA} Let $\MIJ = \MI{\frb ^c}_x \subseteq
\OO$ be \tnparenonly{the germ at $x$ of} any multiplier ideal. If $p
\ge 1$, then no minimal
$p^{\text{th}}$ syzygy   of $\MIJ$ vanishes
modulo 
$\frakm ^{d + 1- p}$. 
\end{theoremalpha}
\noi Let us explain the statement more precisely. 
For the case $p = 1$, fix minimal generators $ f_1, 
\ldots, f_b \in \MIJ$, and let $g_1, \ldots, g_b \in \frakm$ be
functions giving a minimal syzygy
\[  \sum  g_i  f_i \ = \ 0 \] among the $f_i$. Then the claim is that
\[ \ord_x(g_i)
\le d -1\] for at least one index $i$. In general, consider a minimal
free resolution 
\[
\xymatrix{
 \ldots \ar[r]^{u_3} &\OO^{b_2} \ar[r]^{u_2} 
 &   \OO^{b_1} \ar[r]^{u_1} & \OO^{b_0} \ar[r]   & \MIJ \ar[r] & 0
  }
\]
of $\MIJ$, where each $u_p$ is a matrix of elements in $\frakm$ whose
columns minimally generate the module of $p^{\text{th}}$ syzygies  
of
$\MIJ$.  The assertion of the  theorem is that no column
of
$u_p$ (or any
$\CC$-linear combination thereof) can consist entirely of functions
vanishing to order $\ge d + 1 - p$ at $x$. Equivalently, no mimimal
generator of the $p^{\text{th}}$ syzygy module 
\[ \textnormal{Syz}_p(\mathcal{J}) \ =_{\text{def}} \ \Image(u_p) \
\subseteq
\
\OO^{b_{p-1}} \]
of $\mathcal{J}$ lies in $\frakm^{d+ 1 -
p}\cdot \OO^{b_{p-1}}$.

The theorem implies that if $d   \ge 3$, then many integrally closed
ideals cannot arise as multiplier ideals. For  example consider 
$2 \le m \le d-1$   functions \[ f_1, \ldots, f_m \ \in \ \OO \]
vanishing to order
$\ge d$ at
$x$.  If the $f_i$ are chosen generally, then  the complete
intersection ideal
$\mathcal{I} = (f_1, \ldots, f_m)$ that they generate will be radical,
hence integrally closed. On the other hand, the Koszul syzygies among the
$f_i$ violate the condition in Theorem \ref{TheoremA}, and hence
$\mathcal{I}$ is not a multiplier ideal. If $d \ge 3$ a modification
of this construction  (Example \ref {m-primary.example}) yields $\frakm$-primary integrally
closed
 ideals having a syzygy vanishing to high order. The theorem is optimal
when $p = 1$ (Example \ref{Boundary.example}), but by taking into account  all the  $u_i$ we give in \S 2 an extension of Theorem
\ref{TheoremA} that is  generally stronger when $p \ge 2$. 
  Note however that there aren't any restrictions
on the order of vanishing of
\textit{generators} of a multiplier ideal, since for instance all powers
of $\frakm$ occur as multiplier ideals.

Theorem \ref{TheoremA} follows from a more technical statement involving
the vanishing of a map on Koszul cohomology groups. Specifically, 
let $ h_1, \ldots, h_r \in \frakm $ be any
collection of non-zero elements generating an ideal $\fra \subseteq \OO$,
and let
$ K\bull(h_1, \ldots, h_r)$ be the Koszul complex on the $h_i$. We prove:
\begin{theoremalpha} \label{TheoremB}
For
every  $0\le p \le r$, the natural map
\[ H_p\big( K\bull(h_1, \ldots, h_r) \otimes \fra^{r-p}   \MI{\frb^c}
\big)
\lra H_p\big( K\bull(h_1, \ldots, h_r) \otimes  \MI{\frb^c} \big)
\] vanishes.
\end{theoremalpha}
\noi Now fix generators $z_1, \ldots, z_d \in \frakm$, and write $\CC =
\OO/ \frakm$ for the residue field at $x$, viewed as an $\OO$-module.
Taking  $r = d$ and $h_i = z_i$, the theorem implies
\begin{corollaryalpha} \label{CorollaryC} The natural maps
\[
\textnormal{Tor}_p\big(\frakm^{d-p}  \MIJ, \CC\big) \lra
\textnormal{Tor}_p\big( \MIJ, \CC\big )
\] vanish for all $0 \le p \le d$. 
\end{corollaryalpha}
\noi Theorem \ref{TheoremA} is deduced from this statement. 
As for Theorem \ref{TheoremB}, the proof is simply to note that an
exact ``Skoda complex" \cite[Section 9.6.C]{PAG} sits inbetween  the  two Koszul
complexes
 in question. 

We are grateful to J. Lipman, M. Musta\c t\v a and K. Watanabe for
valuable discussions and correspondence. 
 
\section{Proofs}

In this section we prove the results stated in the Introduction. 

 \begin{proof} [Proof of Theorem \ref{TheoremB}] 
Write
$K\bull(h) = K\bull(h_1, \ldots, h_r)$ for the Koszul complex on the
$h_i$, and recall that $\fra = (h_1, \ldots, h_r)$ is the ideal they
generate. This Koszul complex
 contains as a subcomplex the
\textit{Skoda complex}  appearing in
\cite[Section 9.6.C]{PAG}:
\begin{equation}
\xymatrix{
\ldots \ar[r]& \OO^{\binom{r}{2}} \otimes \MI{\fra^{r-2} \frb^c} \ar[r]
&
\OO^r
\otimes
\MI{\fra^{r-1}\frb^c}
\ar[r] &
\MI{\fra^r
\frb^c}
\ar[r] & 0}
. \tag{$ \textnormal{Skod}\bull $}
\end{equation}
 The fact that $ \textnormal{Skod}\bull $ is a subcomplex of $K\bull(h)$
follows from its construction in \cite{PAG}, but more concretely one can
think of
$\textnormal{Skod}\bull$ as being the subcomplex of $K\bull(h)$ arising
by virtue of the inclusions $\fra \cdot \MI{\fra^{\ell -1}\frb^c}
\subseteq
 \MI{\fra^{\ell} \frb^c} $ deduced from equation (*) below.  The basic
fact for our purposes is that
$ \textnormal{Skod}\bull $ is \textit{exact} \cite[Theorem 9.6.36]{PAG}: this
is an elementary consequence of the local vanishing theorems for
multiplier ideals.

Recall next that for any ideals $\frq,\frb \subseteq \OO$, one
has $
\frq \cdot
\MI{\frb^c}
  \subseteq  
\MI{\frq \frb^c}$.  In the case at hand, if we fix $\ell$
this gives inclusions 
\begin{equation} \label{MIinclusions}
\fra^\ell \MI{\frb^c} \
\subseteq \  \MI{\fra ^\ell \frb^c }\
\subseteq \  \MI{\fra ^{\ell-i} \frb^c } \, \subseteq \,  \MI{
\frb^c }  \tag{*}
\end{equation}
for any $i \le \ell$. It follows in the first place that $
\textnormal{Skod}\bull $ is  a subcomplex of $K\bull(h)\otimes
\MI{\frb^c}$. 
Assuming for the time being that $p \ge 1$, let 
\[ \textnormal{Trunc}_p  \ = \ \textnormal{Trunc}_p \big( K\bull(h)
\otimes
\fra^{r-p}
\MI{\frb^c}\big)\
\subseteq
\  K\bull(h)
\otimes
\fra^{r-p}
\MI{\frb^c} \]
be the $p^{\text{th}}$ truncation of $K\bull(h) \otimes \fra^{r-p}
\MI{\frb^p}$, i.e. the complex 
obtained from $K\bull(h) \otimes \fra^{r-p}
\MI{\frb^c}$ by replacing the term $\OO^{\binom{r}{p-1}} \otimes
\fra^{r-p}
\MI{\frb^c}$ by the image $I_{p-1}$ of the incoming map, and
the lower terms by $0$.  Then it also
follows from (*) that
$
\textnormal{Trunc}_p $ is actually a subcomplex of
$\textnormal{Skod}\bull$. Thus all told we have inclusions 
\begin{equation}
\textnormal{Trunc}_p \big( K\bull(h) \otimes
\fra^{r-p}
\MI{\frb^c}\big)  \ \subseteq \ \textnormal{Skod}\bull \ \subseteq \
K\bull(h) \otimes
\MI{\frb^c}, \tag{**} \end{equation}
which are pictured concretely in the following commutative diagram.
\begin{equation*}
\xymatrix@C-10pt{
\ldots \ar[r] & \OO^{\binom{r}{p+1}} \otimes \fra^{r-p} \MI{\frb^c}
\ar[r] \ar@{^{(}->}[d] 
& \OO^{\binom{r}{p}} \otimes \fra^{r-p} \MI{\frb^c}
\ar[r]\ar@{^{(}->}[d]  & I_{p-1}\ar@{^{(}->}[d] 
\ar[r] & \ldots
\\
\ldots \ar[r] & \OO^{\binom{r}{p+1}} \otimes \MI{\fra^{r-p-1} \frb^c}
\ar[r]\ar@{^{(}->}[d]  & \OO^{\binom{r}{p}} \otimes \MI{\fra^{r-p}
\frb^c} \ar[r] \ar@{^{(}->}[d]  &\OO^{\binom{r}{p-1}} \otimes
\MI{\fra^{r-p+1} \frb^c}
\ar[r] \ar@{^{(}->}[d]  &
\ldots\\
\ldots \ar[r] & \OO^{\binom{r}{p+1}} \otimes  \MI{\frb^c}
\ar[r] 
& \OO^{\binom{r}{p}} \otimes  \MI{\frb^c}
\ar[r] & \OO^{\binom{r}{p-1}} \otimes  \MI{\frb^c}
\ar[r] & \ldots } 
\end{equation*}
 
 Theorem \ref{TheoremB} now follows at once from the exactness
of  the Skoda complex. Indeed, 
\[
H_p\big( K\bull(h) \otimes \fra^{r-p}   \MI{\frb^c}
\big) \ = \ H_p \big(
\textnormal{Trunc}_p \big( K\bull(h) \otimes
\fra^{r-p}
\MI{\frb^p}\big)\big)\]
and the map appearing in Theorem \ref{TheoremB} is induced by the
inclusion of the two outer complexes in (**). But this inclusion factors
through the exact complex $ \textnormal{Skod}\bull $, and so the map on
homology vanishes. When $p = 0$ the same argument works ignoring the
truncations (although the statement is tautologous when $p = 0$). 
\end{proof}

\begin{proof}[Proof of Corollary \ref{CorollaryC}] Denote by 
$K\bull \ = \ K\bull(z_1, \ldots, z_d ) $
the Koszul complex on the generators $z_1, \ldots, z_d \in \frakm$. Then
\begin{align*}
\Tor_p\big( \frakm^{d-p} \MI{\frb^c}, \CC\big)  \ &= \ H_p \big(
  K\bull \otimes
\frakm^{d-p}
\MI{\frb^p}\big) \\
\Tor_p\big(  \MI{\frb^c}, \CC\big) \ &= \ H_p \big( K\bull \otimes
\MI{\frb^c} \big),
\end{align*}
and so the assertion  is a special case of Theorem \ref{TheoremB}. 
\end{proof}

Turning to Theorem \ref{TheoremA}, consider any ideal $\II \subseteq \OO$,
and choose a minimal free resolution 
$R\bull$ of $\II$:
\begin{equation} \label{resoln.eqn}
\xymatrix{
 \ldots \ar[r]^{u_3} & R_2 \ar[r]^{u_2} &   R_1 \ar[r]^{u_1} & R_0 
\ar[r]^{\pi}   & \II \ar[r] & 0
  },
\end{equation}
where $R_i = \OO^{b_i}$. 
Fix $p \ge 1$, and let ${e} \in R_p = \OO^{b_p}$ be a generator, so that
$e$ determines a non-zero class in 
\[ \Tor_p(\II, \CC) \ = \ H_p \big( R\bull \otimes \CC\big) \ = \
\CC^{b_p}. \]
In view of Corollary \ref{CorollaryC},  Theorem \ref{TheoremA} follows
from:
\begin{proposition}\label{Tor.Lifting.Prop}
Suppose that there is an integer $a \ge 2$ such that 
\begin{equation} \label{Order.Van.Eqn}
u_p ( e) \ \in \ \frakm^{a}\cdot R_{p-1} \ = \   \frakm^{a} \cdot
\OO^{b_{p-1}}.
\end{equation} Then $e$ represents a class lying in the image of
$\Tor_p(\frakm^{a-1}\II,
\CC)
\lra \Tor_p(\II,
\CC)$. 
\end{proposition}

\begin{proof} We
propose to work explicitly with the identifications
\begin{equation} H_p\big( R\bull \otimes \CC\big) \ = \  \Tor_p(\II, \CC) \
=
\ H_p\big( \II \otimes K\bull
), \notag
\end{equation}
where as above $K\bull = K(z_1, \ldots, z_d)$ is the Koszul complex
\[
\xymatrix{
 \ldots \ar[r]^{\partial_3} & K_2 \ar[r]^{\partial_2} &   K_1
\ar[r]^{\partial_1} & K_0 
\ar[r]   & \CC \ar[r] & 0
  } 
\]
 on
generators $z_1, \ldots , z_d \in \frakm$, so that $K_i =
\OO^{\binom{d}{i}}$.  Specifically, consider the double complex
$R\bull \otimes K\bull$. Then one has isomorphisms
\begin{equation}
 \label{Double.Cx.Isoms}
\begin{aligned}
H_p\big( R\bull \otimes \CC \big) \  &\cong \ H_p \big(
\textnormal{Tot}(R\bull
\otimes K\bull) \big) \\  
H_p \big(
\textnormal{Tot}(R\bull
\otimes K\bull) \big) \ &\cong \ H_p \big(\II \otimes K\bull \big). 
\end{aligned}
\end{equation}
As explained in \cite[Chapter 3.3]{Osborne}, 
the first isomorphism in \eqref{Double.Cx.Isoms} is obtained by using a ``zig-zag"
construction to lift a   generator
$ e   \in    R_p  =   \OO^{b_p} $  to  a
$p$-cycle
\begin{equation} \label{Cyc.in.Total.Cx}
  \widetilde{e} \ = \ \big( e_0, e_1, \ldots, e_p \big) \
\in
\
\textnormal{Tot}(R\bull \otimes K\bull),  \end{equation}
where $e_i \in R_{p-i} \otimes K_i $  and $e_0 =  e \otimes 1 \in   R_p
\otimes K_0$. The second isomorphism in \eqref{Double.Cx.Isoms} arises by
associating to a $p$-cycle $\widetilde{f} = (f_0, \ldots, f_p) \in
\textnormal{Tot}(R\bull
\otimes K\bull)$ the  homology class of the image 
\[ \overline{f}_p \  = \ (\pi \otimes 1)(f_p) \ \in \ \II \otimes K_p\]
of
$f_p
\in R_0
\otimes 
 K_p$.  It suffices to prove:
\begin{quote}
If $u_p(e) \in \frakm^{a} R_{p-1}$, then for $i \ge 1$
one can arrange in \eqref{Cyc.in.Total.Cx} that 
\begin{equation} \label{Ord.Van.Eq.Tot.Cyc}
  e_i \ \in \ \frakm^{a-1} \big( R_{p-i} \otimes K_i \big) 
 \end{equation}
\end{quote}
\noi Indeed, granting this one can take $\overline{e}_p \in
\frakm^{a-1} \II \otimes K_p $. This is a $p$-cycle in $\frakm^{a-1} \II
\otimes K\bull$,  and it  gives the required lifting of
$e$  to a class in 
$\Tor_p(\frakm^{a-1} \II, \CC) = H_p \big( \frakm^{a-1} \II \otimes
K\bull \big)$. 

As for \eqref{Ord.Van.Eq.Tot.Cyc}, consider the commutative diagram
\[
\xymatrix@C+10pt{
\ldots  \ar[r] & R_p \otimes K_2 \ar[d]^{1\otimes \partial_2} \ar[r]^{u_p
\otimes 1} & R_{p-1}
\otimes K_2
\ar[d]^{1\otimes \partial_2}
\ar[r]^{u_{p-1} \otimes 1} & R_{p-2} \otimes K_2 \ar[d]^{1\otimes
\partial_2} \ar[r] &\ldots \\
\ldots \ar[r] & R_p \otimes K_1 \ar[d]^{1\otimes \partial_1} \ar[r]^{u_p
\otimes 1} & R_{p-1}
\otimes K_1 \ar[d]^{1\otimes \partial_1} 
\ar[r]^{u_{p-1}
\otimes 1}  & R_{p-2} \otimes K_1\ar[d]^{1\otimes \partial_1} \ar[r]
&\ldots \\
\ldots \ar[r] & R_p \otimes K_0   \ar[r]^{u_p \otimes 1} & R_{p-1} \otimes
K_0 \ar[r]^{u_{p-1} \otimes 1}  & R_{p-2} \otimes K_0 \ar[r] &\ldots
}
\]
One starts by lifting $
 \pm(u_p\otimes 1)(e \otimes 1) \in R_{p-1} \otimes K_0$
to an element $e_1 \in R_{p-1} \otimes K_1$, i.e. choosing an element
$e_1 \in R_{p-1} \otimes K_1$ such that
\[
  \big( 1 \otimes \partial_1 \big) (e_1) \ = \ (-1)^{p} \cdot (u_p
\otimes 1) (e\otimes 1). \]
The hypothesis on $u_p(e)$ implies that the element on the right lies in
$\frakm^{a}(R_{p-1} \otimes K_0)$. Since $\partial_1$ maps
$\frakm^{a-1} K_1$ onto $\frakm^{a}K_0$, one can
take 
\[ e_1 \ \in \ R_{p-1} \otimes \frakm^{a-1} K_1 \ = \ \frakm^{a-1} \big(
R_{p-1}
\otimes K_1 \big). \] The next step is to lift $\mp( u_{p-1}\otimes
1)(e_1)$ to an element $e_2 \in R_{p-2} \otimes K_2$. 
 The minimality of 
$u_{p-1} : R_{p-1}
\lra R_{p-2}$  implies that 
\[ \big( u_{p-1} \otimes 1\big) (e_1) \ \in \ \frakm^{a} \big( 
R_{p-2}
\otimes K_1 \big) \ = \ \big( R_{p-2} \otimes \frakm^{a} K_{1} \big).
\] Thanks to the exactness of 
$\frakm^{a-1} K_{i+1} \lra \frakm^{a}   K_i \lra \frakm^{a+1}
 K_{i-1}$ for $i \ge 1$, 
 as above we can take
\[  e_2 \ \in \ R_{p-2} \otimes \frakm^{a-1} K_{2}. \]
Continuing in this manner one arrives eventually at
\eqref{Ord.Van.Eq.Tot.Cyc}, completing the proof.  
\end{proof}

\begin{remark} \label{strengthening.order.of.van.prop}
A stronger statement is possible if one takes into account the least
order of vanishing of any of the non-zero entries in a matrix for $u_i$
in the minimal resolution \eqref{resoln.eqn}. Specifically, for
$	1 \le i \le p-1$ let
$\eps_i = \eps_i(\II)$ be the largest integer such that
\begin{equation} 
 \textnormal{Syz}_i(\II) \ = \ \Image(u_i) \ \subseteq \ \frakm^{\eps_i} R_{i-1}. \end{equation}
Then with evident modifications, the argument just completed shows that
in the situation of Proposition \ref{Tor.Lifting.Prop}, 
$ e
\in
\Tor_p(\II, \CC)$ lifts to a class in $\Tor_p(\frakm^{\delta }\II, \CC)$ 
with
\[ \delta \ = \ \delta(\II) \ = \ (a -1) \, + \, (\eps_{p-1} - 1) \, + \,
(\eps_{p-2} - 1) 
\, + \, \ldots \, + \,   (\eps_1 - 1). \qed 
\]
\end{remark}

\section{Examples and Complements}

This section is devoted to some examples and further information
concerning the syzygies of multiplier and other integrally closed ideals.
We begin by outlining a construction of
$\frakm$-primary integrally closed ideals having a minimal first syzygy
that vanishes to high order. The starting point is
\begin{lemma}
Let $J \subseteq \CC[x_1, \dots, x_d]$ be a homogeneous radical ideal,
and let  $\frakm = (x_1, \ldots, x_d)$ be the maximal ideal of
polynomials vanishing at the origin. Then for $k \ge 1$ the ideal
\[  J_k \ =_{\textnormal{def}} \ J \, + \, \frakm^k  \]
is integrally closed.
\end{lemma}
\begin{proof}
Given a homogeneous polynomial $f \in \overline{J_k}$, we need to show
that $f \in J_k$. If $\deg(f) \ge k$ this is trivial, so we can assume
that $f$ has degree $a < k$. By definition $f$ satisfies a polynomial
\[ f^n \, + \, a_1 f^{n-1} \, + \, \ldots \, + \, a_n \ = \ 0, \ \
\text{with} \ \ a_i \in \big ( J + \frakm^k \big)^i. \] We can suppose
that each  $a_i$ is  homogeneous of degree $ai$, and then by considerations
of degree we find that $a_i \in J$. Thus $f^n \in J$, and since $J$ is
radical this implies that $f \in J$. 
\end{proof}

\begin{example}\textbf{(Integrally closed ideal supported at a point with
a syzygy vanishing to high order).} 
\label{m-primary.example}
Keeping notation as in the previous
lemma, let  \[ J \ =\ (f,g)\ \subseteq \ \CC[x_1, \ldots, x_d] \]   be
the complete intersection ideal generated by two general homogeneous
polynomials of   degree $a$. Provided that $d \ge 3$ we can suppose
that $J$ is radical. Now take $k \gg 0$ and set $I = J +
\frakm^k$. Thanks to the lemma $I$ is integrally closed, and hence so too
is its localization 
\[ \mathcal{I}  \ \subseteq \  \OO = \CC[x_1, \ldots,
x_d]_{\frakm}.\] On the other hand, when $k$ is sufficiently
large the Koszul syzygy between $f$ and $g$ remains a minimal first
syzygy of
$\mathcal{I}$. Thus if $d \ge 3$
we have an $\frakm$-primary integrally closed ideal with a  syzygy
vanishing to arbitrary order
$a$ at the origin.  \qed
\end{example}

We next give an example to show that the statement in Theorem
\ref{TheoremA} is optimal when $p = 1$. 
\begin{example} \textbf{(A multiplier ideal on the boundary of Theorem
A).}  \label{Boundary.example} 
Let 
\[  \frb \ = \ (f, g) \ \subseteq \ \CC[x_1, \ldots, x_d] \]
be the complete intersection ideal generated by two general polynomials
vanishing to order $d-1$ at the origin, which we view as defining an ideal
sheaf  on $X = \CC^d$. As above, provided that $d \ge 3$ we can take
$\frb$ to be radical. We claim that 
\begin{equation}\frb \ = \ \MI{X, \frb^2}, \tag{*} \end{equation}
so that we have a multiplier ideal with a first syzygy vanishing to
maximal possible order
$d-1$ at the origin.  To verify (*), let $\mu : X^\pr \lra X$ be the
blow-up of the origin, with exceptional divisor $E$. By choosing $f$ and
$g$ sufficiently generally, we can suppose that
\[ \frb \cdot \OO_{X^\pr} \ = \ \frb^{\pr} \cdot \OO_{X^\pr}(-(d-1)E), \]
where $\frb^\pr \subseteq \OO_{X^\pr}$ is the ideal sheaf of a smooth 
codimension two subvariety meeting $E$ transversely. In particular, 
$\frb^\pr = \MI{X^\pr , (\frb^{\pr})^2 }$. On the other hand, the
birational transformation rule \cite[9.2.33]{PAG} gives that
\begin{align*}
\MI{X, \frb^2} \ &= \ \mu_* \Big( \MMI{X^\pr}{(\frb \cdot
\OO_{X^\pr})^2} \otimes \OO_{X^\pr}(K_{X^\pr/X}) \Big) \\ &= \ \mu_*\Big(
\frb^\pr \cdot \OO_{X^\pr}(-(d-1)E) \Big) \\ &= \   \frb. \qed
\end{align*} 
\end{example}

However when $p \ge 2$, an extension of Theorem \ref{TheoremA} generally
gives a stronger bound. Specifically, combining Corollary \ref{CorollaryC}
with Remark \ref{strengthening.order.of.van.prop}, one arrives at:
\begin{variant} \label{Extension.of.TheoremA}
In the situation of Theorem \ref{TheoremA}, suppose that $\mathcal{J}$
has a minimal $p^{\text{th}}$ syzygy vanishing to order $a_p$ at $x$, and
for $	1 \le i \le p-1$ denote by \[ \eps_i\ = \ \eps_i(\mathcal{J}) \] the least order
of vanishing at $x$ of all non-zero entries in the matrix $u_i$ appearing
in the minimal resolution of $\mathcal{J}$. Then
\begin{equation} 
a_p + \eps_{p-1} + \ldots + \eps_1 \ \le \ d-1. \qed
\tag{*} \end{equation}
\end{variant}

\noi For example, consider when $ d = 4$   the complete intersection ideal
\[  \mathcal{I} \ = \ ( f_1, f_2, f_3) \ \subseteq \ \OO \]
generated by three  functions  vanishing to order $2$ at the
origin. Then $a_2 = \eps_1 =2$, so $\mathcal{I}$ cannot be a multiplier
ideal, but this does not follow from the statement of Theorem
\ref{TheoremA} alone. We do not know how close (*) comes to being optimal.
However the second author constructs in \cite{Leethesis} a number of
examples lying on the boundary of Corollary \ref{CorollaryC}.

Finally, we say a word about adjoint ideals. Let $X$ be a smooth complex
variety, and let 
\[ \fra \ = \ \OO_X(-D) \ \subseteq \OO_X \]
be the ideal sheaf of an integral divisor on $X$. Then $\MI{X, \fra} =
\OO_X(-D)$, so that the multiplier ideal of $D$ is uninteresting. On the
other hand, assuming that $D$ is {reduced}, one can define the
\textit{adjoint ideal}
\[ \textnormal{adj}(D) \ \subseteq \ \OO_X \]
of $D$, which does contain significant information about singularities of
$D$ (see
\cite[Section 9.3.E]{PAG}). One can ask whether these adjont ideals
satisfy the same syzygetic conditions as multiplier ideals. It is a
consequence of the following proposition that this is indeed the case:
\begin{proposition}
Given a reduced divisor $D \subseteq X$, there is an ideal sheaf
$\frb \subseteq \OO_X$ such that
\[ \textnormal{adj}(D) \ = \ \MI{\frb}. \] 
\end{proposition}
\noi One can view this as a converse of \cite[Example 9.3.49]{PAG}. The
Proposition grew out of discussions with Karen Smith and Howard Thompson,
and we thank them for allowing us to include it here. See \cite{SmithThompson} for further applications. 
\begin{proof}[Sketch of Proof]
Let $\mu : X^\pr \lra X$ be a log resolution of $D$, with
\[  \mu^*(D) \ = \ D^\pr + F , \]
$D^\pr$ being the proper transform of $D$. Set $\frb = \mu_* \OO_X(-F)$. We leave it to the reader to check that in fact $\textnormal{adj}(D) = \MI{\frb}$. 
\end{proof}


\begin{thebibliography}{99}

\bibitem{AS}
Urban Angehrn and Yum-Tong Siu, Effective freeness and point separation for adjoint bundles, Invent. Math. \textbf{122} (1995), pp 291--308.


 
\bibitem{Dem1}
Jean-Pierre Demailly, A numerical criterion for very ample line bundles, J. Differential Geom. {\bf 37} (1993), pp.  323--374.

\bibitem{Dem2}
Jean-Pierre Demailly, Multiplier ideal sheaves and analytic methods in algebraic geometry,  in {\it School on Vanishing Theorems and Effective Results in Algebraic Geometry (Trieste, 2000)}, 1--148, Abdus Salam Int. Cent. Theoret. Phys., Trieste, 2001.

\bibitem{EL}
Lawrence  Ein\ and\ Robert Lazarsfeld, A geometric effective Nullstellensatz, Invent. Math. {\bf 137} (1999), pp. 427--448.


\bibitem{ELS}
Lawrence Ein, Robert Lazarsfeld\ and\ Karen E. Smith, Uniform bounds and symmetric powers on a smooth variety, Invent. Math. {\bf 144} (2001), pp.  241--252


\bibitem{FJ}
Charles Favre\ and\ Mattias Jonsson, Valuations and multiplier ideals, J. Amer. Math. Soc. {\bf 18} (2005), pp.  655--684
  
  \bibitem{HM1}
  Christopher Hacon and James McKernan, Boundedness of pluricanonical maps of varieties of general type, preprint 2005.
  
 \bibitem{HM2}
 Christopher Hacon and James McKernan, On the existence of flips, preprint 2005.
  
  \bibitem{PAG}
  Robert Lazarsfeld, 
  \textit{Positivity in Algebraic Geometry, I.-II.} Ergebnisse
der Mathematik und ihrer
  Grenzgebiete, Vols. 48-49., Springer Verlag, Berlin,
2004.

  \bibitem{Leethesis}
  {Kyungyong Lee}, PhD Thesis, University of Michigan, in preparation. 
  
\bibitem{LW}
 Joseph Lipman and Keichi Watanabe, Integrally closed ideals in two-dimensional regular local rings are multiplier ideals,  Math. Res. Lett. {\bf 10} (2003), pp.  423--434;


\bibitem{Osborne}
M. Scott Osborne, {\it Basic Homological Algebra}, Springer, New York, 2000

\bibitem{Siu1} Yum-Tong Siu, Invariance of plurigenera, Invent. Math. \textbf{134} (1998), pp. 661--673.

\bibitem{Siu2}
Yum-Tong Siu, Extension of twisted pluricanonical sections with plurisubharmonic weight and invariance of plurigenera for manifolds not necessarily of general type, in \textit{Complex Geometry, G\"ottingen 2000},  Springer, Berlin (2002), pp. 223--277. 

\bibitem{SmithThompson}
Karen E. Smith and Howard Thompson, Multiplier ideals of schemes admitting a monomial factorizing resolution, in preparation. 

\bibitem{Tak}
Shigeharu Takayama, Pluricanonical systems on algebraic varieties of general type, preprint 2005.


\end{thebibliography}
\end{document}